\documentclass{amsart}

\newtheorem{theorem}{Theorem}[section]
\newtheorem{corollary}[theorem]{Corollary}
\newtheorem{lemma}[theorem]{Lemma}

\begin{document}
\title[Mail your comments to wonniepark@postech.ac.kr]{Normal forms of real hypersurfaces with nondegenerate Levi form}
\author[Won K. Park]{Won K. Park\footnote{ %
E-mail: wonkpark@euclid.postech.ac.kr 
\newline
} \\
..............\\
Any comments, suggestions, errors to\\
wonniepark@postech.ac.kr}
\address{Department of Mathematics, Postech Pohang, Korea, 790-784}
\maketitle

\begin{abstract}
We present a proof of the existence and uniqueness theorem of a normalizing
biholomorphic mapping to Chern-Moser normal form. The explicit form of the
equation of a chain on a real hyperquadric is obtained. There exists a
family of normal forms of real hypersurfaces including Chern-Moser normal
form.
\end{abstract}

\addtocounter{section}{-1}

\section{\textbf{Introduction}}

Let $M$ be an analytic real hypersurface with nondegenerate Levi form in a
complex manifold and $p$ be a point on $M$. Then it is known that there is a
local coordinate system $z^{1},z^{2},\cdots ,z^{n},z^{n+1}\equiv w=u+iv$
with center at $p$, where $M$ is locally defined by the equation 
\begin{equation}
v=\langle z,z\rangle +\sum_{s,t\geq 2}F_{st}(z,\bar{z},u),  \tag*{(0.1)}
\label{0.6}
\end{equation}
where

\begin{enumerate}
\item[(1)]  $\langle z,z\rangle \equiv z^{1}\overline{z^{1}}+\cdots +z^{e}%
\overline{z^{e}}-z^{e+1}\overline{z^{e+1}}-\cdots -z^{n}\overline{z^{n}}$
for a positive integer $e$ in $\frac{n}{2}\leq e\leq n,$

\item[(2)]  $F_{st}(z,\bar{z},u)$ is a real-analytic function of $z,u$ for
each pair $(s,t)\in \Bbb{N}^{2},$ which satisfies 
\begin{equation*}
F_{st}(\mu z,\nu \bar{z},u)=\mu ^{s}\nu ^{t}F_{st}(z,\bar{z},u),
\end{equation*}
for all complex numbers $\mu ,\nu $,

\item[(3)]  the functions $F_{22},F_{23},F_{33}$ satisfy the following
conditions: 
\begin{equation*}
\Delta F_{22}=\Delta ^{2}F_{23}=\Delta ^{3}F_{33}=0,
\end{equation*}
where 
\begin{gather*}
\Delta \equiv D_{1}\overline{D}_{1}+\cdots +D_{e}\overline{D}_{e}-D_{e+1}%
\overline{D}_{e+1}-\cdots -D_{n}\overline{D}_{n}, \\
D_{k}=\frac{\partial }{\partial z^{k}},\quad \overline{D}_{k}=\frac{\partial 
}{\partial \overline{z^{k}}},\quad k=1,\cdots ,n.
\end{gather*}
\end{enumerate}

The local coordinate system \ref{0.6} is called normal coordinate. The
existence of a normal coordinate is a natural consequence of the following
existence theorem of a normalizing biholomorphic mapping to Chern-Moser
normal form.

\begin{theorem}[Chern-Moser]
\label{Thm.4.1}Let $M$ be an analytic real hypersurface with nondegenerate
Levi form at the origin in $\Bbb{C}^{n+1}$ defined by the following
equation: 
\begin{equation*}
v=F(z,\bar{z},u),\quad \left. F\right| _{0}=\left. dF\right| _{0}=0.
\end{equation*}
Then there is a biholomorphic mapping $\phi $ such that 
\begin{equation}
\phi \left( M\right) :v=\langle z,z\rangle +\sum_{s,t\geq 2}F_{st}^{*}(z,%
\bar{z},u)  \tag*{(0.2)}  \label{4.3}
\end{equation}
where 
\begin{equation}
\Delta F_{22}^{*}=\Delta ^{2}F_{23}^{*}=\Delta ^{3}F_{33}^{*}=0. 
\tag*{(0.3)}  \label{4.17}
\end{equation}
\end{theorem}

We shall modify the proof given by Chern and Moser on the existence
theorem(cf. \cite{CM74}) in order that the proof itself yields the
uniqueness theorem as well:

\begin{theorem}[Chern-Moser]
Let $M$ be a nondegenerate analytic real hypersurface in Theorem \ref
{Thm.4.1}. Then the normalization $\phi =(f,g)$ in $\Bbb{C}^{n}\times \Bbb{C}
$ is uniquely determined by the value 
\begin{equation*}
\left. \frac{\partial f}{\partial z}\right| _{0},\quad \left. \frac{\partial
f}{\partial w}\right| _{0},\quad \Re \left( \left. \frac{\partial g}{%
\partial w}\right| _{0}\right) ,\quad \Re \left( \left. \frac{\partial ^{2}g%
}{\partial w^{2}}\right| _{0}\right) .
\end{equation*}
\end{theorem}

An analytic curve on a nondegenerate real hypersurface $M$ is called a chain
if it can be straightened by a normalization. We carefully examine the
equation of a chain on $M$ so that, in particular, we obtain the explicit
form of the equation of a chain $\gamma $ on a real hyperquadric: 
\begin{equation*}
\gamma :\left\{ 
\begin{array}{l}
z=p(\mu ) \\ 
w=\mu +i\langle p(\mu ),p(\mu )\rangle 
\end{array}
\right. 
\end{equation*}
where the function $p(\mu )$ is a solution of the following ordinary
differential equation(cf. \cite{Pi75}) 
\begin{equation*}
p^{\prime \prime }=\frac{2ip^{\prime }\langle p^{\prime },p^{\prime }\rangle
\left( 1+3i\langle p,p^{\prime }\rangle -i\langle p^{\prime },p\rangle
\right) }{\left( 1+i\langle p,p^{\prime }\rangle -i\langle p^{\prime
},p\rangle \right) \left( 1+2i\langle p,p^{\prime }\rangle -2i\langle
p^{\prime },p\rangle \right) }.
\end{equation*}

Normalizations of a real hypersurface $M$ to Chern-Moser normal form is
parameterized by a finite dimensional group $H$ given by 
\begin{equation*}
\left( 
\begin{array}{ccc}
\rho  & 0 & 0 \\ 
-Ca & C & 0 \\ 
-r-i\langle a,a\rangle  & 2ia^{\dagger } & 1
\end{array}
\right) 
\end{equation*}
where 
\begin{equation*}
a^{\dagger }=\left( \overline{a^{1}},\cdots ,\overline{a^{e}},-\overline{%
a^{e+1}},\cdots ,-\overline{a^{n}}\right) .
\end{equation*}
The family of normalization of $M$ shall depend analytically on the
parameters 
\begin{equation*}
C=\left( \left. \frac{\partial f}{\partial z}\right| _{0}\right) ,\quad
-Ca=\left( \left. \frac{\partial f}{\partial w}\right| _{0}\right) ,\quad
\rho =\Re \left( \left. \frac{\partial g}{\partial w}\right| _{0}\right)
,\quad 2\rho r=\Re \left( \left. \frac{\partial ^{2}g}{\partial w^{2}}%
\right| _{0}\right) .
\end{equation*}

We shall show that there is a family of normal forms such that 
\begin{eqnarray*}
v &=&\langle z,z\rangle +\sum_{s,t\geq 2}F_{st}\left( z,\overline{z}%
,u\right) \hspace{1.2in}\text{for }\alpha =0 \\
v &=&-\frac{1}{2\alpha }\ln \left\{ 1-2\alpha \langle z,z\rangle \right\}
+\sum_{s,t\geq 2}F_{st}\left( z,\overline{z},u\right) \quad \text{for }%
\alpha \neq 0
\end{eqnarray*}
where $\alpha \in \Bbb{R}$ and 
\begin{eqnarray*}
\Delta F_{22}\left( z,\overline{z},u\right)  &=&\Delta ^{2}F_{23}\left( z,%
\overline{z},u\right) =0 \\
\Delta ^{3}F_{33}\left( z,\overline{z},u\right)  &=&\beta \Delta ^{4}\left(
F_{22}\left( z,\overline{z},u\right) \right) ^{2}\quad \text{for some }\beta
\in \Bbb{R}.
\end{eqnarray*}
Normalization of a real hypersurface to any normal form among this family is
parametrized by the group $H.$

\section{Existence and Uniqueness Theorem}

Let $M$ be an analytic real hypersurface defined near the origin by 
\begin{equation*}
v=F(z,\bar{z},u),\text{\quad }\left. F\right| _{0}=0 
\end{equation*}
and $\Gamma :\left[ -1,1\right] \rightarrow M$ be an analytic real curve
passing through the origin. Then the equation of $\Gamma $ is given as
follows: 
\begin{equation*}
\Gamma :\left\{ 
\begin{array}{c}
z=p(\mu ), \\ 
w=q(\mu ),
\end{array}
\right. 
\end{equation*}
where $p(0)=q(0)=0.$ Since $\Gamma $ is obviously tangent to $M$ at the
origin, we have the equality 
\begin{equation*}
q^{\prime }(0)=\left( 1+i\left. \frac{\partial F}{\partial u}\right|
_{0}\right) \Re q^{\prime }(0)+i\left( \left. \frac{\partial F}{\partial
z^{\alpha }}\right| _{0}p^{\alpha \prime }(0)+\left. \frac{\partial F}{%
\partial \overline{z}^{\alpha }}\right| _{0}\overline{p}^{\alpha \prime
}(0)\right) . 
\end{equation*}
Then $\Gamma $ is transversal to the complex tangent hyperplane of $M$ at
the origin if and only if the following condition holds 
\begin{equation*}
\Gamma _{*}\frac{d}{dt}\notin \ker \partial \left\{ \frac{w-\overline{w}}{2i}%
-F\left( z,\bar{z},\frac{w+\overline{w}}{2}\right) \right\} . 
\end{equation*}
Then the necessary and sufficient condition for the transversality is given
by the inequality 
\begin{equation*}
\left. \frac{\partial F}{\partial z^{\alpha }}\right| _{0}p^{\alpha \prime
}(0)\neq \frac{1}{2i}\left( 1-i\left. \frac{\partial F}{\partial u}\right|
_{0}\right) q^{\prime }(0). 
\end{equation*}
Thus, under the condition $\left. F_{z}\right| _{0}=\left. F_{\overline{z}%
}\right| _{0}=0,$ the transversality of $\Gamma $ to the complex tangent
hyperplane of $M$ at the origin is equivalent to the inequality 
\begin{equation*}
\Re q^{\prime }(0)=\left\{ 1+\left( \left. \frac{\partial F}{\partial u}%
\right| _{0}\right) ^{2}\right\} ^{-1}q^{\prime }(0)\neq 0. 
\end{equation*}

Hence we suppose that $M$ is an analytic real hypersurface defined by 
\begin{equation}
v=F(z,\bar{z},u),  \tag*{(1.1)}  \label{M}
\end{equation}
where 
\begin{equation*}
\left. F\right| _{0}=\left. F_{z}\right| _{0}=\left. F_{\overline{z}}\right|
_{0}=0. 
\end{equation*}
Let $\Gamma $ be an analytic real curve on $M$ transversal to the complex
tangent hyperplane at the origin. Then the equation of $\Gamma $ is uniquely
given with a distinguished parameterization $\mu $ as follows: 
\begin{equation}
\Gamma :\left\{ 
\begin{array}{l}
z=p(\mu ) \\ 
w=\mu +iF(p(\mu ),\overline{p}(\mu ),\mu )
\end{array}
\right. .  \tag*{(1.2)}  \label{4.34}
\end{equation}

\begin{lemma}
\label{Lem1}Let $g(z,w)$ be a holomorphic function implicitly defined by the
equations: 
\begin{align}
g(z,w)-g(0,w)=& -2iF(p(w),\bar{p}(w),w)  \notag \\
& +2iF\left( z+p(w),\bar{p}(w),w+\frac{1}{2}\{g(z,w)-g(0,w)\}\right) , 
\notag \\
g(0,w)=& iF(p(w),\overline{p}(w),w).  \tag*{(1.3)}  \label{4.60}
\end{align}
Let $\phi $ be a biholomorphic mapping near the origin defined by 
\begin{align}
z=& z^{*}+p(w^{*}),  \notag \\
w=& w^{*}+g(z^{*},w^{*}).  \tag*{(1.4)}  \label{4.61}
\end{align}
Then the mapping $\phi $ transforms the real hypersurface $M$ such that $%
M^{\prime }\equiv \phi \left( M\right) $ is locally defined by an equation
of the following form 
\begin{equation*}
v^{*}=\sum_{\min (s,t)\geq 1}F_{st}^{*}(z^{*},\bar{z}^{*},u^{*})
\end{equation*}
and the curve $\Gamma $ on $M$ via the equation \ref{4.34} is mapped on the $%
u$-curve, $z=v=0.$
\end{lemma}

Note that the holomorphic function $g(z,w)$ is well defined because of the
condition 
\begin{equation*}
\left. F\right| _{0}=\left. F_{z}\right| _{0}=\left. F_{\overline{z}}\right|
_{0}=0, 
\end{equation*}
which implies 
\begin{equation*}
\left. g\right| _{0}=\left. \frac{\partial g}{\partial z}\right| _{0}=\Re
\left( \left. \frac{\partial g}{\partial w}\right| _{0}\right) =0. 
\end{equation*}
Further, the mapping \ref{4.61} is bijective at the origin. Hence the
mapping \ref{4.61} is biholomorphic near the origin for any analytic
function $p(u)$ such that $p(0)=0$.

Suppose that the transformed real hypersurface $M^{\prime }$ is defined by 
\begin{equation*}
v^{*}=F^{*}(z^{*},\overline{z}^{*},u^{*}). 
\end{equation*}
Then the mapping \ref{4.61} yields the following equality: 
\begin{equation}
F(z,\bar{z},u)=F^{*}(z^{*},\bar{z}^{*},u^{*})+\frac{1}{2i}
\{g(z^{*},u^{*}+iv^{*})-\bar{g}(\bar{z}^{*},u^{*}-iv^{*})\},  \tag*{(1.5)}
\label{4.15}
\end{equation}
where 
\begin{eqnarray*}
z & = z^{*}+p(u^{*}+iv^{*}), \\
\bar{z} & = \bar{z}^{*}+\bar{p}(u^{*}-iv^{*}), \\
u & = u^{*}+\frac{1}{2}\{g(z^{*},u^{*}+iv^{*})+\bar{g}(\bar{z}%
^{*},u^{*}-iv^{*})\}.
\end{eqnarray*}
Since $F$ and $F^{*}$ are real analytic, we can consider $z^{*},\bar{z}%
^{*},u^{*}$ as independent variables.\ Hence the condition of $%
F^{*}(z^{*},0,u^{*})=v^{*}=0$ is equivalent via the equality \ref{4.15} to
the following equality: 
\begin{equation}
g(z,u)-\overline{g(0,u)}=2iF\left( z+p(u),\bar{p}(u),u+\frac{1}{2}\{g(z,u)+%
\overline{g(0,u)}\}\right) .  \tag*{(1.6)}  \label{4.50}
\end{equation}
Taking$\ z=0$ yields 
\begin{equation}
g(0,u)-\overline{g(0,u)}=2iF\left( p(u),\bar{p}(u),u+\frac{1}{2}\{g(0,u)+%
\overline{g(0,u)}\}\right) .  \tag*{(1.7)}  \label{4.51}
\end{equation}
Thus we easily see that 
\begin{equation}
g(0,u)+\overline{g(0,u)}=0  \tag*{(1.8)}  \label{4.52}
\end{equation}
if and only if 
\begin{equation*}
g(0,u)=iF(p(u),\bar{p}(u),u). 
\end{equation*}
Let $\Gamma $ be a curve on $M$ defined by a function $p(u)$ via the
equation \ref{4.34}. Then the mapping \ref{4.61} maps the curve $\Gamma $
onto the $u$-curve in $\Bbb{C}^{n+1}$ if and only if the condition \ref{4.52}
on $g(z,w)$ is satisfied.

By requiring the condition \ref{4.52}, the equality \ref{4.50} reduces to 
\begin{eqnarray*}
g(z,u)-g(0,u)= &&-2iF(p(u),\bar{p}(u),u) \\
&&+2iF\left( z+p(u),\bar{p}(u),u+\frac{1}{2}\{g(z,u)-g(0,u)\}\right) .
\end{eqnarray*}
Thus the equalities \ref{4.50} and \ref{4.51} are satisfied by the function $%
g(z,w)$ defined in the equalities \ref{4.60}. This completes the proof of
Lemma \ref{Lem1}.

Note that the mapping \ref{4.61} in Lemma \ref{Lem1} is completely
determined by the analytic function $p(u).$ From the equality \ref{4.60}, we
obtain the expansion of the holomorphic function $g(z,w)$ as a power series
of $z$ up to order 3 inclusive as follows: 
\begin{align}
g(z,w)=& iF(p(w),\overline{p}(w),w)  \notag \\
& +2i(1-iF^{\prime })^{-1}\{F_{\alpha }z^{\alpha }+F_{\alpha \beta
}z^{\alpha }z^{\beta }+F_{\alpha \beta \gamma }z^{\alpha }z^{\beta
}z^{\gamma }\}  \notag \\
& -2(1-iF^{\prime })^{-2}\{F_{\alpha }z^{\alpha }F_{\beta }^{\prime
}z^{\beta }+F_{\alpha }z^{\alpha }F_{\beta \gamma }^{\prime }z^{\beta
}z^{\gamma }  \notag \\
& \hspace{4cm}+F_{\alpha \beta }z^{\alpha }z^{\beta }F_{\gamma }^{\prime
}z^{\gamma }\}  \notag \\
& -2i(1-iF^{\prime })^{-3}\{F_{\alpha }z^{\alpha }(F_{\beta }^{\prime
}z^{\beta })^{2}+2F_{\alpha }z^{\alpha }F_{\beta \gamma }z^{\beta }z^{\gamma
}F^{\prime \prime }  \notag \\
& \hspace{3cm}+(F_{\alpha }z^{\alpha })^{2}F^{\prime \prime }+(F_{\alpha
}z^{\alpha })^{2}F_{\beta }^{\prime \prime }z^{\beta }\}  \notag \\
& +2(1-iF^{\prime })^{-4}\{3(F_{\alpha }z^{\alpha })^{2}F_{\beta }^{\prime
}z^{\beta }F^{\prime \prime }+(F_{\alpha }z^{\alpha })^{3}\cdot F^{\prime
\prime \prime }\}  \notag \\
& +4i(1-iF^{\prime })^{-5}(F_{\alpha }z^{\alpha })^{3}(F^{\prime \prime
})^{2}  \notag \\
& +O(z^{4})  \tag*{(1.9)}  \label{4.59}
\end{align}
where 
\begin{eqnarray*}
F_{\alpha } &=&\sum_{\alpha }\left( \frac{\partial F}{\partial z^{\alpha }}%
\right) \left( p(w),\overline{p}(w),w\right) , \\
F^{\prime } &=&\left( \frac{\partial F}{\partial u}\right) \left( p(w),%
\overline{p}(w),w\right) , \\
F_{\alpha \beta } &=&\frac{1}{2}\sum_{\alpha ,\beta }\left( \frac{\partial
^{2}F}{\partial z^{\alpha }\partial z^{\beta }}\right) \left( p(w),\overline{%
p}(w),w\right) , \\
F_{\alpha }^{\prime } &=&\sum_{\alpha }\left( \frac{\partial ^{2}F}{\partial
z^{\alpha }\partial u}\right) \left( p(w),\overline{p}(w),w\right) , \\
F^{\prime \prime } &=&\frac{1}{2}\left( \frac{\partial ^{2}F}{\partial u^{2}}%
\right) \left( p(w),\overline{p}(w),w\right) , \\
F_{\alpha \beta \gamma } &=&\frac{1}{6}\sum_{\alpha ,\beta ,\gamma }\left( 
\frac{\partial ^{3}F}{\partial z^{\alpha }\partial z^{\beta }\partial
z^{\gamma }}\right) \left( p(w),\overline{p}(w),w\right) , \\
F_{\alpha \beta }^{\prime } &=&\frac{1}{2}\sum_{\alpha ,\beta }\left( \frac{%
\partial ^{3}F}{\partial z^{\alpha }\partial z^{\beta }\partial u}\right)
\left( p(w),\overline{p}(w),w\right) , \\
F_{\alpha }^{\prime \prime } &=&\frac{1}{2}\sum_{\alpha }\left( \frac{%
\partial ^{3}F}{\partial z^{\alpha }\partial u^{2}}\right) \left( p(w),%
\overline{p}(w),w\right) , \\
F^{\prime \prime \prime } &=&\frac{1}{6}\left( \frac{\partial ^{3}F}{%
\partial u^{3}}\right) \left( p(w),\overline{p}(w),w\right) .
\end{eqnarray*}

By Lemma \ref{Lem1}, we have the following condition on the real
hypersurface $M^{\prime }$: 
\begin{equation*}
v=O(z\overline{z}). 
\end{equation*}
Thus it suffices to obtain terms up to $v^{2}$ inclusive in order that we
compute the functions 
\begin{equation*}
F_{st}^{*}(z,\overline{z},u) 
\end{equation*}
of $M^{\prime }$ up to the type $(s,t),$ $s+t\leq 5$ inclusive.

We obtain the expansions of $p^{\alpha }(u+iv)$ and $p^{\bar{\beta}}(u+iv)$
as power series of $v$ as follows: 
\begin{eqnarray*}
p^{\alpha }(u+iv) &=&p^{\alpha }+p^{\alpha \prime }\cdot iv+p^{\alpha \prime
\prime }\cdot \frac{(iv)^{2}}{2}+O(v^{3}) \\
p^{\bar{\beta}}(u+iv) &=&p^{\bar{\beta}}+p^{\bar{\beta}\prime }\cdot iv+p^{%
\bar{\beta}\prime \prime }\cdot \frac{(iv)^{2}}{2}+O(v^{3}).
\end{eqnarray*}
By using this expansion, we expand the holomorphic function $g(z,w)$ as a
power series of $z$ and $v$ in \ref{4.59} as follows: 
\begin{eqnarray*}
g(z,w) &=&\sum_{k,l=0}^{\infty }g_{k}^{(l)}(z,u)\frac{(iv)^{l}}{l!} \\
&=&iF(p(u),\overline{p}(u),u)+g_{0}^{\prime }(0,u)iv+g_{0}^{\prime \prime
}(0,u)\frac{(iv)^{2}}{2} \\
&&+g_{1}(z,u)+g_{1}^{\prime }(z,u)iv+g_{1}^{\prime \prime }(z,u)\frac{%
(iv)^{2}}{2} \\
&&+g_{2}(z,u)+g_{2}^{\prime }(z,u)iv+g_{3}(z,u) \\
&&+O(z^{4})+O(z^{3}v)+O(z^{2}v^{2})+O(zv^{3})+O(v^{3})
\end{eqnarray*}
where 
\begin{equation*}
g_{k}^{(l)}(\mu z,u)=\mu ^{k}g_{k}^{(l)}(z,u)
\end{equation*}
for all complex number $\mu .$

We easily see that the function $g_{k}^{(l)}(z,u)$ depends analytically on
the functions $p(u)$ and $\overline{p}(u),$ polynomially on the derivatives
of $p(u)$ and $\overline{p}(u)$ up to order $l$ inclusive such that the
order sum of the derivatives in each term of $g_{k}^{(l)}(z,u)$ is less than
or equal to the integer $l.$ In low order terms, we obtain 
\begin{eqnarray*}
g_{0}(0,u) &=&iF(p(u),\overline{p}(u),u)=O(u) \\
g_{0}^{\prime }(0,u) &=&iF_{\alpha }p^{\alpha \prime }+iF_{\overline{\beta }%
}p^{\overline{\beta }\prime }+iF^{\prime }=O(1) \\
g_{0}^{\prime \prime }(0,u) &=&iF_{\alpha }p^{\alpha \prime \prime }+iF_{%
\overline{\beta }}p^{\overline{\beta }\prime \prime }+2iF_{\alpha \beta
}p^{\alpha \prime }p^{\beta \prime }+2iF_{\alpha \overline{\beta }}p^{\alpha
\prime }p^{\overline{\beta }\prime }+2iF_{\overline{\alpha }\overline{\beta }%
}p^{\overline{\alpha }\prime }p^{\overline{\beta }\prime } \\
&&+2iF_{\alpha }^{\prime }p^{\alpha \prime }+2iF_{\overline{\beta }}^{\prime
}p^{\overline{\beta }\prime }+2iF^{\prime \prime }
\end{eqnarray*}
and 
\begin{eqnarray*}
g_{1}(z,u) &=&2i(1-F^{\prime })^{-1}\cdot F_{\alpha }z^{\alpha }=O(zu) \\
g_{1}^{\prime }(z,u) &=&2i(1-F^{\prime })^{-1}\left\{ 2F_{\alpha \beta
}z^{\alpha }p^{\beta \prime }+F_{\alpha \overline{\beta }}z^{\alpha }p^{%
\overline{\beta }\prime }+F_{\alpha }^{\prime }z^{\alpha }\right\}  \\
&&-2(1-F^{\prime })^{-1}\left\{ F_{\alpha }^{\prime }p^{\alpha \prime }+F_{%
\overline{\beta }}^{\prime }p^{\overline{\beta }\prime }+2F^{\prime \prime
}\right\} F_{\alpha }z^{\alpha } \\
&=&4i\left. F_{\alpha \beta }\right| _{0}z^{\alpha }p^{\beta \prime
}+2i\left. F_{\alpha \overline{\beta }}\right| _{0}z^{\alpha }p^{\overline{%
\beta }\prime }+O(zu) \\
g_{2}(z,u) &=&2i(1-F^{\prime })^{-1}F_{\alpha \beta }z^{\alpha }z^{\beta
}-2(1-F^{\prime })^{-2}F_{\alpha }z^{\alpha }F_{\beta }^{\prime }z^{\beta }
\\
&&-2i(1-F^{\prime })^{-3}\left( F_{\alpha }z^{\alpha }\right) ^{2}F^{\prime
\prime } \\
&=&2i\left. F_{\alpha \beta }\right| _{0}z^{\alpha }z^{\beta }+O(z^{2}u)
\end{eqnarray*}

The real hypersurface $M^{\prime }$ is defined by the following equation: 
\begin{align}
v=& F\left( z+p(u+iv),\bar{z}+\bar{p}(u-iv),u+\frac{1}{2}\left\{ g(z,u+iv)+%
\bar{g}(\bar{z},u-iv)\right\} \right)   \notag \\
& -\frac{1}{2i}\left\{ g(z,u+iv)-\bar{g}(\bar{z},u-iv)\right\} . 
\tag*{(1.10)}  \label{exp}
\end{align}
We expand the right hand side of the equation \ref{exp} in low order terms
of $v$ as follows: 
\begin{align*}
F& \left( z+p(u+iv),\bar{z}+\bar{p}(u-iv),u+\frac{1}{2}\left\{ g(z,u+iv)+%
\bar{g}(\bar{z},u-iv)\right\} \right)  \\
& \hspace{3.5cm}-\frac{1}{2i}\left\{ g(z,u+iv)-\bar{g}(\bar{z},u-iv)\right\} 
\\
& =A(z,\overline{z},u)+vB(z,\overline{z},u)+v^{2}C(z,\overline{z}%
,u)+O(v^{3}).
\end{align*}
Then we obtain 
\begin{eqnarray*}
A(z,\overline{z},u) &=&F\left( z+p(u),\bar{z}+\bar{p}(u),u+\Re g(z,u)\right)
-\Im g(z,u) \\
B(z,\overline{z},u) &=&iF_{\alpha }\left( z+p(u),\bar{z}+\bar{p}(u),u+\Re
g(z,u)\right) p^{\alpha \prime }(u) \\
&&-iF_{\overline{\beta }}\left( z+p(u),\bar{z}+\bar{p}(u),u+\Re
g(z,u)\right) p^{\overline{\beta }\prime }(u) \\
&&-F^{\prime }\left( z+p(u),\bar{z}+\bar{p}(u),u+\Re g(z,u)\right) \Im
g^{\prime }(z,u) \\
&&-\Re g^{\prime }(z,u) \\
C(z,\overline{z},u) &=&-F_{\alpha \beta }\left( z+p(u),\bar{z}+\bar{p}%
(u),u+\Re g(z,u)\right) p^{\alpha \prime }(u)p^{\beta \prime }(u) \\
&&+F_{\alpha \overline{\beta }}\left( z+p(u),\bar{z}+\bar{p}(u),u+\Re
g(z,u)\right) p^{\alpha \prime }(u)p^{\overline{\beta }\prime }(u) \\
&&-F_{\overline{\alpha }\overline{\beta }}\left( z+p(u),\bar{z}+\bar{p}%
(u),u+\Re g(z,u)\right) p^{\overline{\alpha }\prime }(u)p^{\overline{\beta }%
\prime }(u) \\
&&-iF_{\alpha }^{\prime }\left( z+p(u),\bar{z}+\bar{p}(u),u+\Re
g(z,u)\right) p^{\alpha \prime }(u)\Im g^{\prime }(z,u) \\
&&+iF_{\overline{\beta }}^{\prime }\left( z+p(u),\bar{z}+\bar{p}(u),u+\Re
g(z,u)\right) p^{\overline{\beta }\prime }(u)\Im g^{\prime }(z,u) \\
&&+F^{\prime \prime }\left( z+p(u),\bar{z}+\bar{p}(u),u+\Re g(z,u)\right)
(\Im g^{\prime }(z,u))^{2} \\
&&-\frac{1}{2}F_{\alpha }\left( z+p(u),\bar{z}+\bar{p}(u),u+\Re
g(z,u)\right) p^{\alpha \prime \prime }(u) \\
&&-\frac{1}{2}F_{\overline{\beta }}\left( z+p(u),\bar{z}+\bar{p}(u),u+\Re
g(z,u)\right) p^{\overline{\beta }\prime \prime }(u) \\
&&-\frac{1}{2}F^{\prime }\left( z+p(u),\bar{z}+\bar{p}(u),u+\Re
g(z,u)\right) \Re g^{\prime \prime }(z,u) \\
&&+\frac{1}{2}\Im g^{\prime \prime }(z,u)
\end{eqnarray*}
where 
\begin{eqnarray*}
g^{\prime }(z,u) &=&\left( \frac{\partial g}{\partial w}\right) (z,u) \\
g^{\prime \prime }(z,u) &=&\left( \frac{\partial ^{2}g}{\partial w^{2}}%
\right) (z,u).
\end{eqnarray*}
We decompose the functions $A(z,\overline{z},u),$ $B(z,\overline{z},u),$ $%
C(z,\overline{z},u)$ as follows: 
\begin{eqnarray*}
A(z,\overline{z},u) &=&\sum_{\min (s,t)\geq 1}A_{st}(z,\overline{z},u), \\
B(z,\overline{z},u) &=&\sum_{\min (s,t)\geq 0}B_{st}(z,\overline{z},u), \\
C(z,\overline{z},u) &=&\sum_{\min (s,t)\geq 0}C_{st}(z,\overline{z},u).
\end{eqnarray*}
We easily see the following facts:

\begin{enumerate}
\item[(1)]  $A_{st}(z,\overline{z},u)$ depends analytically on the functions 
$p(u)$ and $\overline{p}(u).$

\item[(2)]  $B_{st}(z,\overline{z},u)$ depends analytically on the function $%
p(u)$ and $\overline{p}(u),$ at most linearly on the derivative $p^{\prime
}(u)$ and $\overline{p}^{\prime }(u).$

\item[(3)]  $C_{st}(z,\overline{z},u)$ depends analytically on the function $%
p(u)$ and $\overline{p}(u)$, at most quadratically on the derivative $%
p^{\prime }(u)$ and $\overline{p}^{\prime }(u),$ at most linearly on the
derivative $p^{\prime \prime }(u)$ and $\overline{p}^{\prime \prime }(u)$
such that the derivative order sum of the derivatives of $p(u)$ and $%
\overline{p}(u)$ in each term is less than or equal to $2.$
\end{enumerate}

\begin{lemma}
\label{c0t}The functions $C_{0t}(z,\overline{z},u),$ $t\in \Bbb{N},$ do not
depend on the derivative $p^{\prime \prime }(u)$ and $\overline{p}^{\prime
\prime }(u).$
\end{lemma}

This claim is easily verified by observing that the following functions 
\begin{equation*}
C(z,\overline{z},u)\quad \text{and\quad }-\frac{1}{2}\left( \frac{\partial
^{2}A}{\partial u^{2}}\right) (z,\overline{z},u) 
\end{equation*}
depend in the same manner on the derivative $p^{\prime \prime }(u)$ and $%
\overline{p}^{\prime \prime }(u),$ because 
\begin{equation*}
A(z,\overline{z},u+iv)=A(z,\overline{z},u)+iv\left( \frac{\partial A}{%
\partial u}\right) (z,\overline{z},u)-\frac{v^{2}}{2}\left( \frac{\partial
^{2}A}{\partial u^{2}}\right) (z,\overline{z},u)+\cdots . 
\end{equation*}
By the way, $A(z,0,u)=0$ is the defining equation of the function $g(z,u)$
in Lemma \ref{Lem1} with $g(0,u)=iF(p(u),\overline{p}(u),u).$ Thus we have
the following identities 
\begin{equation*}
\left( \frac{\partial A}{\partial u}\right) (z,0,u)=\left( \frac{\partial
^{2}A}{\partial u^{2}}\right) (z,0,u)=\cdots =0. 
\end{equation*}
Note that the identity 
\begin{equation*}
\left( \frac{\partial ^{2}A}{\partial u^{2}}\right) (z,0,u)=0 
\end{equation*}
gives the desired relation between the terms having $p^{\prime \prime }$ and 
$\overline{p}^{\prime \prime }$ and the terms not having $p^{\prime \prime }$
and $\overline{p}^{\prime \prime }$ so that we verify the function $C(z,0,u)$
is independent of the derivative $p^{\prime \prime }(u)$ and $\overline{p}%
^{\prime \prime }(u).$ This proves the claim in Lemma \ref{c0t}.

Explicitly, we compute the expansion of the right hand side of the equation 
\ref{exp} in low order terms so that 
\begin{align*}
v=& A_{11}(z,\overline{z},u)+A_{22}(z,\overline{z},u)+A_{12}(z,\overline{z}%
,u)+A_{13}(z,\overline{z},u)+A_{23}(z,\overline{z},u) \\
& \text{\hspace{1.5cm}}+A_{21}(z,\overline{z},u)+A_{31}(z,\overline{z}%
,u)+A_{32}(z,\overline{z},u)\} \\
& +v\{B_{00}(z,\overline{z},u)+B_{11}(z,\overline{z},u)+B_{01}(z,\overline{z}%
,u)+B_{02}(z,\overline{z},u)+B_{12}(z,\overline{z},u) \\
& \hspace{1.5cm}+B_{10}(z,\overline{z},u)+B_{20}(z,\overline{z},u)+B_{21}(z,%
\overline{z},u)\} \\
& +v^{2}\{C_{00}(z,\overline{z},u)+C_{01}(z,\overline{z},u)+C_{10}(z,%
\overline{z},u)\} \\
& +O\left( z^{1}\overline{z}^{4}\right) +O\left( z^{4}\overline{z}%
^{1}\right) +O\left( vz^{3}\right) +O\left( v\overline{z}^{3}\right)
+O\left( v^{3}\right)  \\
& +\sum_{\min (s,t)\geq 1,s+t\geq 6}O\left( z^{s}\overline{z}^{t}\right)
+\sum_{s+t\geq 4}O\left( vz^{s}\overline{z}^{t}\right) +\sum_{s+t\geq
2}O\left( v^{2}z^{s}\overline{z}^{t}\right) ,
\end{align*}
where 
\begin{eqnarray*}
A_{11}(z,\overline{z},u) &=&F_{\alpha \overline{\beta }}z^{\alpha }z^{%
\overline{\beta }}-i\left( 1+iF^{\prime }\right) ^{-1}F_{\alpha }^{\prime
}z^{\alpha }F_{\overline{\beta }}z^{\overline{\beta }} \\
&&+i\left( 1-iF^{\prime }\right) ^{-1}F_{\alpha }z^{\alpha }F_{\overline{%
\beta }}^{\prime }z^{\overline{\beta }} \\
&&+2\left( 1+iF^{\prime }\right) ^{-1}\left( 1-iF^{\prime }\right)
^{-1}F^{\prime \prime }F_{\alpha }z^{\alpha }F_{\overline{\beta }}z^{%
\overline{\beta }} \\
&=&\langle z,z\rangle +O(z\overline{z}u) \\
B_{00}(z,\overline{z},u) &=&i\left( 1+iF^{\prime }\right) ^{-1}F_{\alpha
}p^{\alpha \prime }-i\left( 1-iF^{\prime }\right) ^{-1}F_{\overline{\beta }%
}p^{\overline{\beta }\prime }-\left( F^{\prime }\right) ^{2} \\
&=&O(1) \\
B_{01}(z,\overline{z},u) &=&2iF_{\alpha \overline{\beta }}p^{\alpha \prime
}z^{\overline{\beta }}+2F^{\prime \prime }F_{\overline{\alpha }}z^{\overline{%
\alpha }}+i\left( 1+iF^{\prime }\right) F_{\overline{\alpha }}z^{\overline{%
\alpha }} \\
&&+iF_{\alpha }^{\prime }p^{\alpha \prime }\left( 1+iF^{\prime }\right)
^{-1}F_{\overline{\alpha }}z^{\overline{\alpha }} \\
&&+2i\left( F_{\alpha }p^{\alpha \prime }+F_{\overline{\beta }}^{\prime }p^{%
\overline{\beta }\prime }\right) F^{\prime \prime }\left( 1+iF^{\prime
}\right) ^{-1}F_{\overline{\gamma }}z^{\overline{\gamma }}
\end{eqnarray*}

Then we obtain 
\begin{align*}
v& =F_{11}^{*}(z,\overline{z},u)+F_{22}^{*}(z,\overline{z},u)+F_{12}^{*}(z,%
\overline{z},u)+F_{13}^{*}(z,\overline{z},u)+F_{23}^{*}(z,\overline{z},u) \\
& \text{\hspace{1.5cm}}+F_{21}^{*}(z,\overline{z},u)+F_{31}^{*}(z,\overline{z%
},u)+F_{32}^{*}(z,\overline{z},u) \\
& \text{\hspace{1.5cm}}+O(z\overline{z}^{4})+O(z^{4}\overline{z})+\sum_{\min
(s,t)\geq 1,s+t\geq 6}O(z^{s}\overline{z}^{t}),
\end{align*}
where 
\begin{eqnarray*}
F_{11}^{*} &=&(1-B_{00})^{-1}A_{11} \\
F_{12}^{*} &=&(1-B_{00})^{-1}A_{12}+(1-B_{00})^{-2}A_{11}B_{01} \\
F_{13}^{*}
&=&(1-B_{00})^{-1}A_{13}+(1-B_{00})^{-2}(A_{11}B_{02}+A_{12}B_{01}) \\
&&+(1-B_{00})^{-3}A_{11}B_{01}^{2} \\
F_{22}^{*}
&=&(1-B_{00})^{-1}A_{22}+(1-B_{00})^{-2}(A_{11}B_{11}+A_{12}B_{10}+A_{21}B_{01})
\\
&&+(1-B_{00})^{-3}(2A_{11}B_{01}B_{10}+A_{11}^{2}C_{00}) \\
F_{23}^{*}
&=&(1-B_{00})^{-1}A_{23}+(1-B_{00})^{-2}(A_{11}B_{12}+A_{12}B_{11}+A_{21}B_{02}
\\
&&\hspace{6cm}+A_{13}B_{10}+A_{22}B_{01}) \\
&&+(1-B_{00})^{-3}(2A_{11}B_{01}B_{11}+2A_{11}B_{10}B_{02}+2A_{12}B_{01}B_{10}
\\
&&\hspace{3cm}+A_{21}B_{01}^{2}+A_{11}^{2}C_{01}+2A_{11}A_{12}C_{00}) \\
&&+3(1-B_{00})^{-4}(A_{11}B_{01}^{2}B_{10}+A_{11}^{2}B_{01}C_{00}).
\end{eqnarray*}

By Lemma \ref{c0t}, the functions $F_{22}^{*},F_{23}^{*}$ does not depend on
the derivative $p^{\prime \prime }$ and $\overline{p}^{\prime \prime },$ and
the dependence of the coefficients in $F_{22}^{*},F_{23}^{*}$ on the
derivative $p^{\prime }$ and $\overline{p}^{\prime }$ is of the form: 
\begin{equation}
\frac{A_{1}(u,p,\overline{p},p^{\prime },\overline{p}^{\prime })}{%
(1-B_{00})^{3}},  \tag*{(1.11)}  \label{4.71}
\end{equation}
and 
\begin{equation}
\frac{A_{2}(u,p,\overline{p},p^{\prime },\overline{p}^{\prime })}{%
(1-B_{00})^{4}},  \tag*{(1.12)}  \label{4.72}
\end{equation}
where $A_{1}$ depends analytically on $u,p,\overline{p}$ and at most
quadratically on $p^{\prime },\overline{p}^{\prime }$ and $A_{2}$ depends
analytically on $u,p,\overline{p},$ at most cubically on $p^{\prime },%
\overline{p}^{\prime }$.

For future reference, we analyze the terms containing the first order
derivatives $p^{\prime },\overline{p}^{\prime }$ in $B_{00}$ and $B_{01}$ so
that 
\begin{align*}
O(up^{\prime })+O(u\overline{p}^{\prime })\hspace{2cm}\text{in}& \quad
B_{00}(z,\overline{z},u) \\
2i\langle p^{\prime },z\rangle +O(\overline{z}up^{\prime })+O(\overline{z}u%
\overline{p}^{\prime })\quad \text{in}& \quad B_{01}(z,\overline{z},u).
\end{align*}
Thus analyzing the terms containing the second order derivatives $p^{\prime
\prime }$ and $\overline{p}^{\prime \prime }$ in 
\begin{equation*}
\left( \frac{\partial F_{11}^{*}}{\partial u}\right) (z,\overline{z},u)\quad 
\text{and}\quad \left( \frac{\partial F_{12}^{*}}{\partial u}\right) (z,%
\overline{z},u), 
\end{equation*}
we obtain 
\begin{align}
O(z\overline{z}up^{\prime \prime })+O(z\overline{z}u\overline{p}^{\prime
\prime })\hspace{2cm}\text{in}& \quad \left( \frac{\partial F_{11}^{*}}{%
\partial u}\right) (z,\overline{z},u)  \notag \\
2i\langle z,z\rangle \langle p^{\prime \prime },z\rangle +O(z\overline{z}%
^{2}up^{\prime \prime })+O(z\overline{z}^{2}u\overline{p}^{\prime \prime
})\quad \text{in}& \quad \left( \frac{\partial F_{12}^{*}}{\partial u}%
\right) (z,\overline{z},u).  \tag*{(1.13)}  \label{4.95}
\end{align}

Then we are ready to present a proof of the existence theorem for
Chern-Moser normal form.

\begin{theorem}[Chern-Moser]
\label{CM}There is a biholomorphic mapping $\phi $ which transforms $M$ to a
real hypersurface of the following form: 
\begin{equation}
v=\langle z,z\rangle +\sum_{\min (s,t)\geq 2}F_{st}(z,\bar{z},u),  \notag
\end{equation}
where 
\begin{equation}
\Delta ^{2}F_{23}=0.  \tag*{(1.14)}  \label{f23}
\end{equation}
Geometrically, there exists a unique analytic curve $\Gamma $ on $M$ which
passes through the origin and is tangent to a vector transversal to the
complex tangent hyperplane at the origin and which is mapped onto the $u$%
-curve by the biholomorphic mapping $\phi .$ Further, there exists a
biholomorphic mapping $\phi _{1}$ which, in addition to \ref{f23}, achieves
the following conditions: 
\begin{equation*}
\Delta F_{22}=\Delta ^{3}F_{33}=0.
\end{equation*}
\end{theorem}

Let $M^{\prime }$ be a real hypersurface obtained in Lemma \ref{Lem1} by the
biholomorphic mapping \ref{4.61}, which is defined by the following
equation: 
\begin{equation*}
v=\sum_{\min (s,t)\geq 1}F_{st}^{*}(z,\bar{z},u). 
\end{equation*}
Then there is a unique analytic function $D(z,u)$(cf. \cite{CM74}) such that 
\begin{equation*}
F_{11}^{*}\left( z+D(z,u),\overline{z},u\right) =\sum_{s\geq 1}F_{s1}^{*}(z,%
\bar{z},u), 
\end{equation*}
and the function $D(z,u)$ satisfies the condition 
\begin{equation*}
D(0,u)=D_{z}(0,u)=0. 
\end{equation*}
Thus $D(z,u)$ depends analytically of $u,p,\overline{p}$ and rationally of
the derivative $p^{\prime },\overline{p}^{\prime }.$

We decompose the function $D(z,u)$ such that 
\begin{equation*}
D(z,u)=\sum_{s\geq 2}D_{s}(z,u),
\end{equation*}
where 
\begin{equation*}
D_{s}(\mu z,u)=\mu ^{s}D_{s}(z,u)\text{\quad for all }\mu \in \Bbb{C}.
\end{equation*}
Then the functions $D_{2}(z,u),D_{3}(z,u)$ are given by 
\begin{eqnarray*}
A_{11}\left( D_{2}(z,u),\overline{z},u\right) 
&=&A_{21}+(1-B_{00})^{-1}A_{11}B_{10} \\
A_{11}\left( D_{3}(z,u),\overline{z},u\right) 
&=&A_{31}+(1-B_{00})^{-1}(A_{11}B_{20}+A_{21}B_{10}) \\
&&+(1-B_{00})^{-2}A_{11}B_{10}^{2}.
\end{eqnarray*}
Note that $D_{2}(z,u),D_{3}(z,u)$ do not depend of the second order
derivative $p^{\prime \prime },\overline{p}^{\prime \prime }.$

Then we obtain 
\begin{eqnarray*}
v &=&\sum_{\min (s,t)\geq 1}F_{st}^{*}\left( z,\bar{z},u\right)  \\
&=&F_{11}^{*}\left( z,\overline{z},u\right) +F_{11}^{*}\left( z,\overline{%
D(z,u)},u\right) +F_{11}^{*}\left( D(z,u),\overline{z},u\right)  \\
&&+\sum_{\min (s,t)\geq 2}F_{st}^{*}\left( z,\bar{z},u\right)  \\
&=&F_{11}^{*}\left( z+D(z,u),\overline{z+D(z,u)},u\right) +\sum_{\min
(s,t)\geq 2}G_{st}\left( z,\bar{z},u\right) .
\end{eqnarray*}
We notice 
\begin{eqnarray*}
G_{22}\left( z,\bar{z},u\right)  &=&F_{22}^{*}\left( z,\bar{z},u\right)
-F_{11}^{*}\left( D_{2}(z,u),\overline{D_{2}(z,u)},u\right)  \\
G_{23}\left( z,\bar{z},u\right)  &=&F_{23}^{*}\left( z,\bar{z},u\right)
-F_{11}^{*}\left( D_{2}(z,u),\overline{D_{3}(z,u)},u\right) .
\end{eqnarray*}
We easily see that the functions $G_{22},G_{23}$ depend on $u,p,\overline{p}%
,p^{\prime },\overline{p}^{\prime }$ in the same form as respectively in \ref
{4.71} and \ref{4.72}.

We take an analytic function $E(u)$ such that 
\begin{equation*}
F_{11}^{*}(z,\overline{z},u)=\langle E(u)z,E(u)z\rangle ,\quad \text{and}%
\quad E(0)=id_{n\times n}. 
\end{equation*}
Note that the function $E(u)$ is determined up to the following relation: 
\begin{equation*}
E_{1}(u)=U(u)E(u), 
\end{equation*}
where 
\begin{equation}
\langle U(u)z,U(u)z\rangle =\langle z,z\rangle ,\quad \text{and}\quad
U(0)=id_{n\times n}.  \tag*{(1.15)}  \label{4.40}
\end{equation}
Then the biholomorphic mapping defined by the following equation: 
\begin{eqnarray*}
z^{*} &=&E(w)\{z+D(z,w)\}, \\
w^{*} &=&w,
\end{eqnarray*}
transforms $M^{\prime }$ to a real hypersurface of the following form: 
\begin{equation}
v=\langle z,z\rangle +\sum_{\min (s,t)\geq 2}H_{st}(z,\bar{z},u). 
\tag*{(1.16)}  \label{4.5}
\end{equation}
By the way, we still obtain a real hypersurface in \ref{4.5} by a
biholomorphic mapping as follows: 
\begin{align}
z^{*}& =U(w)E(w)\{z+D(z,w)\},  \notag \\
w^{*}& =w,  \tag*{(1.17)}  \label{4.6}
\end{align}
where the holomorphic function $U(w)$ satisfy the condition \ref{4.40}.

By using the expansion 
\begin{eqnarray*}
E(u) &=&E(w)-ivE^{\prime }(w)+\cdots  \\
U(u) &=&U(w)-ivU^{\prime }(w)+\cdots ,
\end{eqnarray*}
we obtain 
\begin{align}
v=& F_{11}^{*}\left( z+D(z,u),\overline{z+D(z,u)},u\right) +\sum_{\min
(s,t)\geq 2}G_{st}(z,\bar{z},u)  \notag \\
=& \left\langle E(u)(z+D(z,u)),E(u)(z+D(z,u))\right\rangle +\sum_{\min
(s,t)\geq 2}G_{st}(z,\bar{z},u)  \notag \\
=& \left\langle U(u)E(u)(z+D(z,u)),U(u)E(u)(z+D(z,u))\right\rangle
+\sum_{\min (s,t)\geq 2}G_{st}(z,\bar{z},u)  \notag \\
=& \left\langle U(w)E(w)(z+D(z,w)),U(w)E(w)(z+D(z,w))\right\rangle   \notag
\\
& -iv\left\langle U^{\prime
}(w)E(w)(z+D(z,w)),U(w)E(w)(z+D(z,w))\right\rangle   \notag \\
& +iv\left\langle U(w)E(w)(z+D(z,w)),U^{\prime
}(w)E(w)(z+D(z,w))\right\rangle   \notag \\
& -iv\left\langle \{E^{\prime
}(w)(z+D(z,w))+E(w)D_{u}(z,w)\},E(w)(z+D(z,w))\right\rangle   \notag \\
& +iv\left\langle E(w)(z+D(z,w)),\{E^{\prime
}(w)(z+D(z,w))+E(w)D_{u}(z,w)\}\right\rangle   \notag \\
& +O(z\overline{z}v^{2})+\sum_{\min (s,t)\geq 2}G_{st}(z,\bar{z},u), 
\tag*{(1.18)}  \label{toe}
\end{align}
where 
\begin{gather*}
w=u+iv, \\
U^{\prime }(u)=\frac{dU}{du}(u),\quad E^{\prime }(u)=\frac{dE}{du}(u), \\
D_{u}(z,w)=\left( \frac{\partial D}{\partial u}\right) (z,w).
\end{gather*}
By introducing a holomorphic variable $z^{\natural }=z+D(z,w),$ we obtain
from the equation \ref{toe}: 
\begin{align}
v=& \left\langle U(w)E(w)z^{\natural },U(w)E(w)z^{\natural }\right\rangle
+G_{22}(z^{\natural },\bar{z}^{\natural },u)  \notag \\
& -i\left\langle E(u)z^{\natural },E(u)z^{\natural }\rangle \{\langle
U^{\prime }(u)E(u)z^{\natural },U(u)E(u)z^{\natural }\right\rangle   \notag
\\
& -\left\langle U(u)E(u)z^{\natural },U^{\prime }(u)E(u)z^{\natural
}\right\rangle   \notag \\
& -i\left\langle E(u)z^{\natural },E(u)z^{\natural }\right\rangle \left\{
\left\langle E^{\prime }(u)z^{\natural },E(u)z^{\natural }\right\rangle
-\left\langle E(u)z^{\natural },E^{\prime }(u)z^{\natural }\right\rangle
\right\}   \notag \\
& +G_{23}^{*}(z^{\natural },\bar{z}^{\natural },u).+G_{32}^{*}(z^{\natural },%
\bar{z}^{\natural },u).  \notag \\
& +\sum_{\min (s,t)\geq 2,s+t\geq 6}G_{st}^{*}(z^{\natural },\bar{z}%
^{\natural },u).  \tag*{(1.19)}  \label{4.96}
\end{align}
where 
\begin{align}
G_{23}^{*}(z,\bar{z},u)=& G_{23}(z,\bar{z},u)+iF_{11}^{*}(z,\overline{z}%
,u)F_{11}^{*}\left( z,\overline{\left( \frac{\partial D_{2}}{\partial u}%
\right) (z,u)},u\right)   \notag \\
& -\sum_{\beta }\left( \frac{\partial G_{22}}{\partial \overline{z}^{\beta }}%
\right) (z,\bar{z},u)\overline{D_{2}^{\beta }(z,u)}  \notag \\
=& G_{23}(z,\bar{z},u)+iF_{11}^{*}(z,\overline{z},u)\left( \frac{\partial
F_{12}^{*}}{\partial u}\right) \left( z,\overline{z},u\right)   \notag \\
& -iF_{11}^{*}(z,\overline{z},u)\left( \frac{\partial F_{11}^{*}}{\partial u}%
\right) \left( z,\overline{D_{2}(z,u)},u\right)   \notag \\
& -\sum_{\beta }\left( \frac{\partial G_{22}}{\partial \overline{z}^{\beta }}%
\right) (z,\bar{z},u)\overline{D_{2}^{\beta }(z,u)}  \tag*{(1.20)}
\label{4.12}
\end{align}

By the equalities \ref{4.95}, the dependence of the functions $G_{23}^{*}$
on the second order derivatives $p^{\prime \prime }$ and $\overline{p}%
^{\prime \prime }$ is given as follows 
\begin{equation*}
-2\langle z,z\rangle ^{2}\langle p^{\prime \prime },z\rangle +O(z^{2}%
\overline{z}^{3}up^{\prime \prime })+O(z^{2}\overline{z}^{3}u\overline{p}%
^{\prime \prime })\quad \text{in}\quad G_{23}^{*}(z,\overline{z},u). 
\end{equation*}
Notice that the function $G_{23}^{*}(z,\bar{z},u)$ in \ref{4.96} and \ref
{4.12} is independent of the function $U(u).$

Therefore after the biholomorphic mapping in \ref{4.6}, we obtain 
\begin{eqnarray*}
v= &&\langle z^{*},z^{*}\rangle +H_{22}(z^{*},\bar{z}^{*},u)+H_{23}(z^{*},%
\bar{z}^{*},u).+H_{32}(z^{*},\bar{z}^{*},u). \\
&&+\sum_{\min (s,t)\geq 2,s+t\geq 6}H_{st}(z^{*},\bar{z}^{*},u).
\end{eqnarray*}
where 
\begin{eqnarray*}
H_{22}(z,\bar{z},u)= &&G_{22}\left( E(u)^{-1}U(u)^{-1}z,\overline{%
E(u)^{-1}U(u)^{-1}z},u\right)  \\
&&-i\langle z,z\rangle \{\langle U^{\prime }(u)U(u)^{-1}z,z\rangle -\langle
z,U^{\prime }(u)U(u)^{-1}z\rangle \} \\
&&-i\langle z,z\rangle \{\langle E^{\prime
}(u)E(u)^{-1}U(u)^{-1}z,U(u)^{-1}z\rangle  \\
&&\hspace{2cm}-\langle U(u)^{-1}z,E^{\prime }(u)E(u)^{-1}U(u)^{-1}z\rangle \}
\end{eqnarray*}
and the dependence of $H_{23}(z,\bar{z},u)$ on $p^{\prime \prime },\overline{%
p}^{\prime \prime }$ is as follows: 
\begin{equation*}
H_{23}(z,\bar{z},0)=-2\langle z,z\rangle ^{2}\langle p^{\prime \prime
}(0),z\rangle +K_{23}(z,\overline{z},0;p^{\prime }(0),\overline{p^{\prime
}(0)}).
\end{equation*}

By using the following identity 
\begin{equation*}
\Delta ^{2}\{\langle z,z\rangle ^{2}\langle p,z\rangle \}=2(n+1)(n+2)\langle
p,z\rangle , 
\end{equation*}
the equation $\Delta ^{2}H_{23}=0$ is a second order ordinary differential
equation 
\begin{equation*}
A_{1}p^{\prime \prime }+A_{2}\overline{p}^{\prime \prime }=B 
\end{equation*}
where

\begin{enumerate}
\item[(1)]  $A_{1},A_{2}$ are $n\times n$ matrix valued functions and $B$ is 
$\Bbb{C}^{n}$-valued function,

\item[(2)]  $A_{1}=id_{n\times n}+O(u)$ and $A_{2}=O(u),$

\item[(3)]  $A_{1},A_{2},B$ depend analytically of $u,p,\overline{p},$

\item[(4)]  $A_{1},A_{2}$ depend at most linearly of $p^{\prime },\overline{p%
}^{\prime },$

\item[(5)]  $B$ depends at most cubically of $p^{\prime },\overline{p}%
^{\prime }.$
\end{enumerate}

\noindent Then we obtain 
\begin{align}
p^{\prime \prime }& =Q(u,p,\overline{p},p^{\prime },\overline{p}^{\prime }) 
\notag \\
& \equiv \left( A_{1}-A_{2}\overline{A_{1}^{-1}}\overline{A_{2}}\right)
^{-1}\left( B-A_{2}\overline{A_{1}^{-1}}\overline{B}\right)  \tag*{(1.21)}
\label{puleaf-eq}
\end{align}
where the function $Q$ depends rationally on the derivatives $p^{\prime },%
\overline{p}^{\prime }.$

Therefore there exists a unique analytic curve $\Gamma $ on $M$ which passes
through the origin and is tangent to a vector transversal to the complex
tangent hyperplane at the origin and which is mapped by a biholomorphic
mapping into the $u$-curve.

Since $\left\langle U(u)z,U(u)z\right\rangle =\langle z,z\rangle ,$ we have
identities 
\begin{gather*}
\left\langle U^{\prime }(u)U(u)^{-1}z,z\right\rangle +\left\langle
z,U^{\prime }(u)U(u)^{-1}z\right\rangle =0 \\
\text{\textrm{Tr}}(U^{\prime }(u)U(u)^{-1})+\overline{\text{\textrm{Tr}}%
(U^{\prime }(u)U(u)^{-1})}=0,
\end{gather*}
where 
\begin{equation*}
\text{\textrm{Tr}}(A)=\mathrm{trace}\text{\textrm{\ of\quad }}z\rightarrow
Az,
\end{equation*}
Then the equation $\Delta H_{22}=0$ is given as follows: 
\begin{align}
& \left\langle U(u)^{-1}U^{\prime }(u)z,z\right\rangle +\frac{1}{2(n+2)}%
\langle z,z\rangle \text{\textrm{Tr}}(U(u)^{-1}U^{\prime }(u))  \notag \\
=& \frac{1}{2i(n+2)}\Delta G_{22}\left( E(u)^{-1}z,\overline{E(u)^{-1}z}%
,u\right)   \notag \\
& -\frac{1}{2}\left\{ \left\langle E^{\prime }(u)E(u)^{-1}z,z\right\rangle
-\left\langle z,E^{\prime }(u)E(u)^{-1}z\right\rangle \right\}   \notag \\
& -\frac{1}{2(n+2)}\langle z,z\rangle \left\{ \text{\textrm{Tr}}(E^{\prime
}(u)E(u)^{-1})-\overline{\text{\textrm{Tr}}(E^{\prime }(u)E(u)^{-1})}%
\right\} .  \tag*{(1.22)}  \label{4.94}
\end{align}
By using the following identities 
\begin{gather*}
\Delta \{\langle z,z\rangle \langle Az,z\rangle \}=(n+2)\langle Az,z\rangle +%
\text{\textrm{Tr}}(A)\langle z,z\rangle  \\
\Delta ^{2}\{\langle z,z\rangle \langle Az,z\rangle \}=2(n+1)\text{\textrm{Tr%
}}(A),
\end{gather*}
we obtain 
\begin{eqnarray*}
\langle z,z\rangle \text{\textrm{Tr}}(U(u)^{-1}U^{\prime }(u))= &&\frac{1}{%
4i(n+1)}\Delta ^{2}G_{22}\left( E(u)^{-1}z,\overline{E(u)^{-1}z},u\right)  \\
&&-\frac{1}{2}\left\{ \text{\textrm{Tr}}(E^{\prime }(u)E(u)^{-1})-\overline{%
\text{\textrm{Tr}}(E^{\prime }(u)E(u)^{-1})}\right\} .
\end{eqnarray*}
Thus the equation \ref{4.94} is a first order ordinary differential equation
of $U(u)$ as follows: 
\begin{eqnarray*}
&&\left\langle U(u)^{-1}U^{\prime }(u)z,z\right\rangle  \\
= &&\frac{1}{2i(n+2)}\Delta G_{22}\left( E(u)^{-1}z,\overline{E(u)^{-1}z}%
,u\right)  \\
&&-\frac{1}{8i(n+1)(n+2)}\langle z,z\rangle \Delta ^{2}G_{22}\left(
E(u)^{-1}z,\overline{E(u)^{-1}z},u\right)  \\
&&-\frac{1}{2}\left\{ \left\langle E^{\prime }(u)E(u)^{-1}z,z\right\rangle
-\left\langle z,E^{\prime }(u)E(u)^{-1}z\right\rangle \right\}  \\
&&-\frac{1}{4(n+2)}\langle z,z\rangle \left\{ \text{\textrm{Tr}}(E^{\prime
}(u)E(u)^{-1})-\overline{\text{\textrm{Tr}}(E^{\prime }(u)E(u)^{-1})}%
\right\} .
\end{eqnarray*}
Hence by requiring 
\begin{equation*}
U(0)=E(0)=id_{n\times n},
\end{equation*}
there is a unique biholomorphic mapping 
\begin{align}
z^{*}& =U(w)E(w)\{z+D(z,w)\},  \notag \\
w^{*}& =w,  \tag*{(1.23)}  \label{4.33}
\end{align}
which transforms $M^{\prime }$ to a real hypersurface of the following form: 
\begin{equation}
v=\langle z,z\rangle +\sum_{\min (s,t)\geq 2}H_{st}(z,\bar{z},u) 
\tag*{(1.24)}  \label{4.36}
\end{equation}
where 
\begin{equation*}
\Delta H_{22}=\Delta H_{23}=0.
\end{equation*}

We consider the following mappings 
\begin{align}
\phi _{1}& :\left\{ 
\begin{array}{l}
z=z^{*}+p(w^{*}) \\ 
w=w^{*}+g(z^{*},w^{*})
\end{array}
\right.  \notag \\
\phi _{2}& :\left\{ 
\begin{array}{l}
z^{*}=E(w)(z+D(z,w)) \\ 
w^{*}=w
\end{array}
\right.  \tag*{(1.25)}  \label{p123} \\
\phi _{3}& :\left\{ 
\begin{array}{l}
z^{*}=\sqrt{\mathrm{sign}\{q^{\prime }(0)\}q^{\prime }(w)}Uz \\ 
w^{*}=q(w)
\end{array}
\right.  \notag
\end{align}
where $p(w),$ $g(z,w),$ $E(w),$ $D(z,w),$ $q(w)$ are holomorphic functions
satisfying 
\begin{gather*}
\overline{g(0,u)}=-g(0,u),\quad \overline{q(u)}=q(u), \\
p(0)=q(0)=0,\quad \det q^{\prime }(0)\neq 0,\quad \det U\neq 0 \\
E(0)=id_{n\times n},\quad D(0,w)=D_{z}(0,w)=0.
\end{gather*}
We easily see by parameter counting that the mapping 
\begin{equation}
(\phi _{1},\phi _{2},\phi _{3})\longmapsto \phi _{3}\circ \phi _{2}\circ
\phi _{1}  \tag*{(1.26)}  \label{pi123}
\end{equation}
is bijective. Hence a biholomorphic mapping $\phi ,$ $\left. \phi \right|
_{0}=0,$ has a unique decomposition 
\begin{equation*}
\phi =\phi _{3}\circ \phi _{2}\circ \phi _{1}. 
\end{equation*}

Note that $U(w)=E(w)=id_{n\times n}$ and $D(z,w)=0$ in \ref{4.33} whenever $%
M $ is already in the form \ref{4.36}. Hence, from the decomposition \ref
{pi123}, we easily see that any biholomorphic mapping, preserving the form 
\ref{4.36} and the $u$-curve, is given by 
\begin{align}
z^{*}& =\sqrt{\text{\textrm{sign}}\{q^{\prime }(0)\}q^{\prime }(w)}Uz, 
\notag \\
w^{*}& =q(w),  \tag*{(1.27)}  \label{4.93}
\end{align}
where 
\begin{gather*}
\overline{q(w)}=q(\bar{w}),\text{\quad }q(0)=0,\text{\quad }q^{\prime
}(0)\neq 0, \\
U\in GL(n;\Bbb{C}),\text{\quad }\langle Uz,Uz\rangle =\text{\textrm{sign}}%
\{q^{\prime }(0)\}\langle z,z\rangle .
\end{gather*}

The mapping in \ref{4.93} transforms the real hypersurface defined by 
\begin{eqnarray*}
v^{*}= &&\langle z^{*},z^{*}\rangle +H_{22}^{*}(z^{*},\bar{z}%
^{*},u^{*})+H_{23}^{*}(z^{*},\bar{z}^{*},u^{*})+H_{32}^{*}(z^{*},\bar{z}%
^{*},u^{*}) \\
&&+H_{33}^{*}(z^{*},\bar{z}^{*},u^{*})+O(z^{*4}\bar{z}^{*2})+O(z^{*2}\bar{z}%
^{*4}) \\
&&+\sum_{\min (s,t)\geq 2,s+t\geq 7}H_{st}^{*}(z^{*},\bar{z}^{*},u^{*})
\end{eqnarray*}
to a real hypersurface as follows: 
\begin{eqnarray*}
v= &&\langle z,z\rangle +q^{\prime }H_{22}^{*}(Uz,\overline{Uz},q(u)) \\
&&+q^{\prime }\sqrt{\left| q^{\prime }\right| }\{H_{23}^{*}(Uz,\overline{Uz}%
,q(u))+H_{32}^{*}(Uz,\overline{Uz},q(u))\} \\
&&+sq^{\prime 2}H_{33}^{*}(Uz,\overline{Uz},q(u))+\left\{ \frac{1}{2}\left( 
\frac{q^{\prime \prime }}{q^{\prime }}\right) ^{2}-\frac{q^{\prime \prime
\prime }}{3q^{\prime }}\right\} \langle z,z\rangle ^{3} \\
&&+O(z^{4}\bar{z}^{2})+O(z^{2}\bar{z}^{4}) \\
= &&\langle z,z\rangle +H_{22}(z,\overline{z},u)+H_{23}(z,\overline{z}%
,u)+H_{32}(z,\overline{z},u)+H_{33}(z,\overline{z},u) \\
&&+O(z^{4}\bar{z}^{2})+O(z^{2}\bar{z}^{4}).
\end{eqnarray*}
Hence we obtain 
\begin{align}
H_{22}(z,\overline{z},u)& =q^{\prime }H_{22}^{*}(Uz,\overline{Uz},q(u)) 
\notag \\
H_{23}(z,\overline{z},u)& =q^{\prime }\sqrt{\left| q^{\prime }\right| }%
H_{23}^{*}(Uz,\overline{Uz},q(u))  \notag \\
H_{33}(z,\overline{z},u)& =sq^{\prime 2}H_{33}^{*}(Uz,\overline{Uz}%
,q(u))+\left\{ \frac{1}{2}\left( \frac{q^{\prime \prime }}{q^{\prime }}%
\right) ^{2}-\frac{q^{\prime \prime \prime }}{3q^{\prime }}\right\} \langle
z,z\rangle ^{3}.  \tag*{(1.28)}  \label{4.92}
\end{align}
Note that $\Delta H_{22}^{*}=\Delta ^{2}H_{23}^{*}=0$ whenever $\Delta
H_{22}=\Delta ^{2}H_{23}=0.$

We can achieve the condition $\Delta ^{3}H_{33}^{*}=0$ by a third order
ordinary differential equation as follows: 
\begin{equation}
\frac{q^{\prime \prime \prime }}{3q^{\prime }}-\frac{1}{2}\cdot \left( \frac{%
q^{\prime \prime }}{q^{\prime }}\right) ^{2}=\kappa (u),  \tag*{(1.29)}
\label{4.43}
\end{equation}
where 
\begin{equation*}
\kappa (u)=-\frac{1}{6n(n+1)(n+2)}\cdot \Delta ^{3}H_{33}(z,\overline{z},u). 
\end{equation*}
The differential equation in \ref{4.43} determines a projective parameter on
the $u$-curve. This completes the proof of Theorem \ref{CM}.

\begin{theorem}[Chern-Moser]
\label{CMU}Let $M$ be a nondegenerate analytic real hypersurface defined by
the equation 
\begin{equation*}
v=F\left( z,\overline{z},u\right) \quad \left. F\right| _{0}=\left.
dF\right| _{0}=0.
\end{equation*}
Then a biholomorphic normalizing mapping of $M,$ $\phi =(f,g)$ in $\Bbb{C}%
^{n}\times \Bbb{C},$ is uniquely determined by the value 
\begin{equation*}
\left. \frac{\partial f}{\partial z}\right| _{0},\quad \left. \frac{\partial
f}{\partial w}\right| _{0},\quad \Re \left( \left. \frac{\partial g}{%
\partial w}\right| _{0}\right) ,\quad \Re \left( \left. \frac{\partial ^{2}g%
}{\partial w^{2}}\right| _{0}\right) .
\end{equation*}
\end{theorem}

As noted above, a biholomorphic mapping $\phi $ satisfying $\left. \phi
\right| _{0}=0$ is uniquely decomposed to 
\begin{equation*}
\phi =\phi _{3}\circ \phi _{2}\circ \phi _{1} 
\end{equation*}
where $\phi _{1},\phi _{2},\phi _{3}$ are biholomorphic mappings in \ref
{p123} satisfying 
\begin{equation*}
\left. \phi _{1}\right| _{0}=\left. \phi _{2}\right| _{0}=\left. \phi
_{3}\right| _{0}=0. 
\end{equation*}
Note that the mapping $(\phi _{1},\phi _{2},\phi _{3})\mapsto \phi =\phi
_{3}\circ \phi _{2}\circ \phi _{1}$ is bijective.

We take $\phi $ to be a normalizing biholomorphic mapping of $M.$ Then the
uniqueness of $\phi _{1},$ $\phi _{2},$ $\phi _{3}$ up to the value 
\begin{equation*}
\left. \frac{\partial f}{\partial z}\right| _{0},\quad \left. \frac{\partial
f}{\partial w}\right| _{0},\quad \Re \left( \left. \frac{\partial g}{%
\partial w}\right| _{0}\right) ,\quad \Re \left( \left. \frac{\partial ^{2}g%
}{\partial w^{2}}\right| _{0}\right) 
\end{equation*}
assures the uniqueness of the normalizing mapping $\phi $. The uniqueness of 
$\phi _{1},$ $\phi _{2},$ $\phi _{3}$ each is verified in the proof of
Theorem \ref{CM} through uniquely determining the holomorphic functions 
\begin{equation*}
p(w),\quad E(w),\quad q(w) 
\end{equation*}
via the ordinary differential equations \ref{puleaf-eq}, \ref{4.94}, \ref
{4.43} by the initial values 
\begin{equation*}
p^{\prime }(0),\quad E(0)\equiv U,\quad q^{\prime }(0),\quad q^{\prime
\prime }(0). 
\end{equation*}

We may have the following relations: 
\begin{eqnarray*}
\sqrt{\left| q^{\prime }(0)\right| }U &=\left. \frac{\partial f}{\partial z}%
\right| _{0} \\
-\sqrt{\left| q^{\prime }(0)\right| }Up^{\prime }(0) &=\left( 1-i\left. 
\frac{\partial F}{\partial u}\right| _{0}\right) ^{-1}\left. \frac{\partial f%
}{\partial w}\right| _{0} \\
q^{\prime }(0) &=\Re \left( \left. \frac{\partial g}{\partial w}\right|
_{0}\right) \\
2q^{\prime }(0)q^{\prime \prime }(0) &=\Re \left\{ \left( 1-i\left. \frac{%
\partial F}{\partial u}\right| _{0}\right) ^{-2}\left. \frac{\partial ^{2}g}{%
\partial w^{2}}\right| _{0}\right\} .
\end{eqnarray*}
For the case $\left. dF\right| _{0}=0$ rather than $\left. F_{z}\right|
_{0}=\left. F_{\overline{z}}\right| _{0}=0,$ we have simpler relations: 
\begin{eqnarray*}
\sqrt{\left| q^{\prime }(0)\right| }U &=\left. \frac{\partial f}{\partial z}%
\right| _{0},\quad -\sqrt{\left| q^{\prime }(0)\right| }Up^{\prime
}(0)=\left. \frac{\partial f}{\partial w}\right| _{0} \\
q^{\prime }(0) &=\Re \left( \left. \frac{\partial g}{\partial w}\right|
_{0}\right) ,\quad 2q^{\prime }(0)q^{\prime \prime }(0)=\Re \left( \left. 
\frac{\partial ^{2}g}{\partial w^{2}}\right| _{0}\right)
\end{eqnarray*}
so that the values $p^{\prime }(0),$ $E(0)\equiv U,$ $q^{\prime }(0),$ $%
q^{\prime \prime }(0)$ are uniquely determined by 
\begin{equation*}
\left. \frac{\partial f}{\partial z}\right| _{0},\quad \left. \frac{\partial
f}{\partial w}\right| _{0},\quad \Re \left( \left. \frac{\partial g}{%
\partial w}\right| _{0}\right) ,\quad \Re \left( \left. \frac{\partial ^{2}g%
}{\partial w^{2}}\right| _{0}\right) . 
\end{equation*}
This completes the proof of Theorem \ref{CMU}.

\section{Chains and Orbit Parameters}

\textbf{I}. Let $M$ be a nondegenerate analytic real hypersurface. Then we
may define a family of distinguished curves on $M$ via Chern-Moser normal
form, which are defined alternatively and identified to be the same by E.
Cartan \cite{Ca32} and Chern-Moser \cite{CM74}. Let $\gamma
:(0,1)\rightarrow M$ be an open connected curve. Then the curve $\gamma $ is
called a chain if, for each point $p\in \gamma ,$ there exist an open
neighborhood $U$ of the point $p$ and a biholomorphic mapping $\phi $ on $U$
which translates the point $p$ to the origin and transforms $M$ to
Chern-Moser normal form such that 
\begin{equation*}
\phi \left( U\cap \gamma \right) \subset \left\{ z=v=0\right\} . 
\end{equation*}
By Theorem \ref{CM} and Theorem \ref{CMU}, a chain $\gamma $ locally exists
uniquely for each vector transversal to the complex tangent plane such that $%
\gamma $ is tangential to the vector.

From the proof of Theorem \ref{CM}, we have an ordinary differential
equation which locally characterizes a chain $\gamma ,$ passing through the
origin $0\in M$. Suppose that 
\begin{equation*}
\gamma :\left\{ 
\begin{array}{l}
z=p(u) \\ 
w=u+iF\left( p(u),\overline{p}(u),u\right)
\end{array}
\right. 
\end{equation*}
Then there exists an ordinary differential equation 
\begin{equation}
p^{\prime \prime }=Q\left( u,p,\overline{p},p^{\prime },\overline{p}^{\prime
}\right)  \tag*{(2.1)}  \label{ordinary}
\end{equation}
such that the function $p(u)$ is a solution of the ordinary differential
equation \ref{ordinary}.

We take $M$ to be the real hyperquadric $v=\langle z,z\rangle .$ Then the
chain $\gamma $ is locally given by 
\begin{equation*}
\gamma :\left\{ 
\begin{array}{l}
z=p(u) \\ 
w=u+i\langle p(u),p(u)\rangle
\end{array}
\right. . 
\end{equation*}
With $F\left( z,\overline{z},u\right) =\langle z,z\rangle ,$ we obtain the
equation $\Delta ^{2}F_{23}=0$ as follows 
\begin{equation*}
\left\{ \left( 1-i\langle p^{\prime },p\rangle +i\langle p,p^{\prime
}\rangle \right) p^{\prime \prime }-ip^{\prime }\langle p^{\prime \prime
},p\rangle \right\} +ip^{\prime }\langle p,p^{\prime \prime }\rangle
=2ip^{\prime }\langle p^{\prime },p^{\prime }\rangle . 
\end{equation*}
Then we easily check that the equation of chains on a real hyperquadric is
given by 
\begin{equation*}
p^{\prime \prime }=\frac{2ip^{\prime }\langle p^{\prime },p^{\prime }\rangle
\left( 1+3i\langle p,p^{\prime }\rangle -i\langle p^{\prime },p\rangle
\right) }{\left( 1+i\langle p,p^{\prime }\rangle -i\langle p^{\prime
},p\rangle \right) \left( 1+2i\langle p,p^{\prime }\rangle -2i\langle
p^{\prime },p\rangle \right) }. 
\end{equation*}

\textbf{II}. The isotropy subgroup of the automorphism group of a real
hyperquadric $v=\langle z,z\rangle $ consists of fractional linear mappings $%
\phi $ such that 
\begin{equation}
\phi =\phi _{\sigma }:\left\{ 
\begin{array}{c}
z^{*}=\frac{C(z-aw)}{1+2i\langle z,a\rangle -w(r+i\langle a,a\rangle )} \\ 
w^{*}=\frac{\rho w}{1+2i\langle z,a\rangle -w(r+i\langle a,a\rangle )}
\end{array}
\right.  \tag*{(2.2)}  \label{2.1}
\end{equation}
where the constants $\sigma =(C,a,\rho ,r)$ satisfy 
\begin{gather*}
a\in \Bbb{C}^{n},\quad \rho \neq 0,\quad \rho ,r\in \Bbb{R}, \\
C\in GL(n;\Bbb{C}),\quad \langle Cz,Cz\rangle =\rho \langle z,z\rangle .
\end{gather*}
Further, $\phi $ decomposes to 
\begin{equation*}
\phi =\varphi \circ \psi , 
\end{equation*}
where 
\begin{equation}
\psi :\left\{ 
\begin{array}{c}
z^{*}=\frac{z-aw}{1+2i\langle z,a\rangle -i\langle a,a\rangle w} \\ 
w^{*}=\frac{w}{1+2i\langle z,a\rangle -i\langle a,a\rangle w}
\end{array}
\right. \quad \text{and}\quad \varphi :\left\{ 
\begin{array}{c}
z^{*}=\frac{Cz}{1-rw} \\ 
w^{*}=\frac{\rho w}{1-rw}
\end{array}
\right. .  \tag*{(2.3)}  \label{2.2}
\end{equation}
Hence the local automorphisms of a real hyperquadric is identified with a
group $H$ of the following matrices: 
\begin{equation*}
\left( 
\begin{array}{ccc}
\rho & 0 & 0 \\ 
-Ca & C & 0 \\ 
-r-i\langle a,a\rangle & 2ia^{\dagger } & 1
\end{array}
\right) 
\end{equation*}
where 
\begin{equation*}
a^{\dagger }=\left( \overline{a^{1}},\cdots ,\overline{a^{e}},-\overline{%
a^{e+1}},\cdots ,-\overline{a^{n}}\right) . 
\end{equation*}

We easily verify 
\begin{equation*}
\phi _{\sigma }^{*}(v-\langle z,z\rangle )=(v-\langle z,z\rangle )\rho
(1+\delta )^{-1}(1+\bar{\delta})^{-1}, 
\end{equation*}
where 
\begin{equation*}
1+\delta =1+2i\langle z,a\rangle -(r+i\langle a,a\rangle )w. 
\end{equation*}
Hence the automorphisms $\phi _{\sigma }$ are normalizations of a real
hyperquadric. Further, by Theorem \ref{CMU}, each normalization of a real
hyperquadric is necessarily an automorphism. Then a chain $\gamma $ on a
real hyperquadric is necessarily given by 
\begin{eqnarray*}
\gamma &=\phi ^{-1}\left( z=v=0\right) \\
&=\left\{ \left( \frac{\rho ^{-1}ua}{1-\rho ^{-1}u\left( -r+i\langle
a,a\rangle \right) },\frac{\rho ^{-1}u}{1-\rho ^{-1}u\left( -r+i\langle
a,a\rangle \right) }\right) \right\} \\
&=\left\{ v=\langle z,z\rangle \right\} \cap \Bbb{C}\left( a,1\right)
\end{eqnarray*}
so that the chain $\gamma $ is just an intersection of a complex line.

By Theorem \ref{CMU}, each normalization $N=(f,g)$ is uniquely determined by
the initial value 
\begin{equation*}
C,\quad a,\quad \rho ,\quad r 
\end{equation*}
such that 
\begin{eqnarray*}
f(z,w) &=C(z-aw)+f^{*}(z,w) \\
g(z,w) &=\rho (w+rw^{2})+g^{*}(z,w)
\end{eqnarray*}
where 
\begin{equation*}
\left. f^{*}\right| _{0}=\left. df^{*}\right| _{0}=\left. g^{*}\right|
_{0}=\left. dg^{*}\right| _{0}=\Re \left( \left. g_{ww}^{*}\right|
_{0}\right) =0. 
\end{equation*}
Hence the group $H=\left\{ \left( C,a,\rho ,r\right) \right\} $
parameterizes the normalizations of a real hypersurface. Further, Theorem 
\ref{CM} and Theorem \ref{CMU} together yields a family of polynomial
identities(cf. \cite{Pa2}). Then we have showed that the group $H$ gives a
group action via normalization on the class of normalized real hypersurfaces.

\textbf{III}. Suppose that $M$ is an analytic real hypersurface defined near
the origin by the following equation: 
\begin{equation*}
v=\langle z,z\rangle +\sum_{\alpha ,\beta }\left( \kappa _{\alpha \beta
}z^{\alpha }z^{\beta }+\kappa _{\overline{\alpha }\overline{\beta }}z^{%
\overline{\alpha }}z^{\overline{\beta }}\right) +F(z,\bar{z},u) 
\end{equation*}
where 
\begin{equation*}
F(z,\bar{z},u)=\sum_{k=3}^{\infty }F_{k}\left( z,z,u\right) . 
\end{equation*}
Let $N_{\sigma }$ be a normalization of $M$ with the initial value $\sigma
=(C,a,\rho ,r)\in H$ and let $\phi _{\sigma ^{\prime }}=\varphi \circ \psi $
be a local automorphism of a real hyperquadric(cf. \ref{2.1}, \ref{2.2})
with the initial value $\sigma ^{\prime }=(C,a,\rho ,r_{0})\in H,$ where 
\begin{equation*}
r_{0}=r-\Re (\kappa _{\alpha \beta }a^{\alpha }a^{\beta }). 
\end{equation*}
Then there are two decompositions of $N_{\sigma }$ as follows(cf. \cite{CM74}%
): 
\begin{equation*}
\left\{ 
\begin{array}{l}
N_{\sigma }=E\circ \phi _{\sigma ^{\prime }}=E\circ \varphi \circ \psi , \\ 
N_{\sigma }=\varphi \circ E\circ \psi ,
\end{array}
\right. 
\end{equation*}
where $E$ is the normalization with the identity initial value$.$

Let $M$ be a nondegenerate real hypersurface and $N_{\sigma }$ be a
normalization of $M$ with initial value $\sigma =(C,a,\rho ,r)\in H$ such
that $M^{\prime }\equiv N_{\sigma }\left( M\right) $ is defined by the
equation 
\begin{equation*}
v=\langle z,z\rangle +\sum_{\min (s,t)\geq 2}F_{st}^{*}(z,\bar{z},u) 
\end{equation*}
where 
\begin{equation*}
\Delta F_{22}^{*}=\Delta ^{2}F_{23}^{*}=\Delta ^{3}F_{33}^{*}=0. 
\end{equation*}
Note that the mapping $\varphi $ is itself a normalization in the following
decomposition: 
\begin{align}
N_{\sigma }& =\varphi \circ E\circ \psi  \notag \\
& =\varphi \circ N_{\sigma _{1}}  \tag*{(2.4)}  \label{2.30}
\end{align}
where $N_{\sigma _{1}}$ is a normalization of $M$ with initial value $\sigma
_{1}=(id_{n\times n},a,1,0)\in H.$

As a consequences of the decomposition \ref{2.30}, we notice that a
normalization $N_{\sigma }$ is analytic of 
\begin{equation*}
z,\quad w,\quad C,\quad \rho ,\quad r 
\end{equation*}
near the point $z=w=r=0$ and $C=id_{n\times n},\rho =1.$ More precisely, 
\begin{equation*}
N_{\sigma }=\left( \frac{Cf(z,w)}{1-rg(z,w)},\frac{\rho g(z,w)}{1-rg(z,w)}%
\right) 
\end{equation*}
where 
\begin{equation*}
N_{\sigma _{1}}=\left( f(z,w),g(z,w)\right) \quad \sigma _{1}=(id_{n\times
n},a,1,0). 
\end{equation*}
Notice that the size of convergence of the normalization $N_{\sigma }$ at
the origin is determined by the value $a,r.$

Further, suppose that the transformed real hypersurface $N_{\sigma }\left(
M\right) $ is defined by 
\begin{equation*}
v=\langle z,z\rangle +F^{*}\left( z,\bar{z},u\right) . 
\end{equation*}
Then the function $F^{*}\left( z,\bar{z},u\right) $ is real-analytic of 
\begin{equation*}
z,\quad w,\quad C,\quad \rho ,\quad r 
\end{equation*}
near the point $z=w=r=0$ and $C=id_{n\times n},\rho =1.$

In fact, from the proof of Theorem \ref{CM}, we obtain

\begin{theorem}
Let $M$ be a nondegenerate real-analytic real hypersurface defined by 
\begin{equation*}
v=F\left( z,\overline{z},u\right) \quad \left. F\right| _{0}=\left.
F_{z}\right| _{0}=\left. F_{\overline{z}}\right| _{0}=0.
\end{equation*}
Then $N_{\sigma }$ and $F^{*}\left( z,\overline{z},u\right) $ are analytic
of 
\begin{equation*}
z,\quad w,\quad C,\quad a,\quad \rho ,\quad r.
\end{equation*}
\end{theorem}

We have a natural group action by normalizations on the class of real
hypersurfaces in normal form(cf. \cite{Pa2}). Then, under a natural
compact-open topology(cf. \cite{Na}), we obtain

\begin{corollary}
The local automorphism group of a nondegenerate analytic real hypersurface
is a Lie group.
\end{corollary}

\section{Normal forms of real hypersurfaces}

\textbf{I}. Let $M$ be a nondegenerate analytic real hypersurface defined by
the equation 
\begin{equation*}
v=F_{11}^{*}\left( z,\overline{z},u\right) +\sum_{s,t\geq 2}F_{st}^{*}\left(
z,\overline{z},u\right) .
\end{equation*}
Let's take a matrix valued function $E(u)$ satisfying 
\begin{equation*}
F_{11}^{*}\left( z,\overline{z},u\right) =\langle E(u)z,E(u)z\rangle .
\end{equation*}
Then suppose that the biholomorphic mapping 
\begin{equation*}
\phi :\left\{ 
\begin{array}{l}
z^{*}=E\left( w\right) z \\ 
w^{*}=w
\end{array}
\right. 
\end{equation*}
transforms $M$ to another real hypersurface $\phi \left( M\right) $ defined
by 
\begin{equation*}
v=F_{11}^{*}\left( z,\overline{z},u\right) +\sum_{s,t\geq 2}G_{st}\left( z,%
\overline{z},u\right) .
\end{equation*}
Then we have the following relation 
\begin{eqnarray*}
&&F_{11}^{*}\left( E(u)^{-1}E^{\prime }(u)z,\overline{z},u\right)  \\
= &&-\frac{2i}{n+2}\cdot \mathrm{tr}G_{22}\left( z,\overline{z},u\right) +%
\frac{i}{(n+1)(n+2)}\cdot (\mathrm{tr})^{2}G_{22}\left( z,\overline{z}%
,u\right) \cdot F_{11}^{*}\left( z,\overline{z},u\right)  \\
&&+\frac{2i}{n+2}\cdot \mathrm{tr}F_{22}^{*}\left( E_{1}(u)^{-1}E(u)z,%
\overline{E_{1}(u)^{-1}E(u)z},u\right)  \\
&&-\frac{i}{(n+1)(n+2)}\cdot (\mathrm{tr})^{2}F_{22}^{*}\left( z,\overline{z}%
,u\right) \cdot F_{11}^{*}\left( z,\overline{z},u\right)  \\
&&+\frac{1}{2}\left( \frac{\partial F_{11}^{*}}{\partial u}\right) \left( z,%
\overline{z},u\right) .
\end{eqnarray*}
where $E_{1}(u)$ is a given matrix valued function satisfying 
\begin{equation*}
F_{11}^{*}\left( z,\overline{z},u\right) =\langle E_{1}(u)z,E_{1}(u)z\rangle
.
\end{equation*}

Here we have constant solutions 
\begin{equation*}
E(u)=E(0) 
\end{equation*}
whenever 
\begin{equation*}
\mathrm{tr}F_{22}^{*}\left( z,\overline{z},u\right) =\mathrm{const}%
.F_{11}^{*}\left( z,\overline{z},u\right) =\mathrm{tr}G_{22}\left( z,%
\overline{z},u\right) 
\end{equation*}
and 
\begin{equation*}
\left( \frac{\partial F_{11}^{*}}{\partial u}\right) \left( z,\overline{z}%
,u\right) =0. 
\end{equation*}
As a necessary condition for a normal form, we require that the $u$-curve be
a chain. From the observation in the previous paragraph, a normal form may
have to take the following form 
\begin{equation*}
v=\langle z,z\rangle +\sum_{s,t\geq 2}F_{st}^{*}\left( z,\overline{z}%
,u\right) 
\end{equation*}
where 
\begin{eqnarray*}
\Delta F_{22}^{*}\left( z,\overline{z},u\right) &=\mathrm{const.}\langle
z,z\rangle \\
\Delta ^{2}F_{23}^{*}\left( z,\overline{z},u\right) &=0.
\end{eqnarray*}
Suppose that a real hypersurface $M$ is defined by the equation 
\begin{align}
v& =\langle z,z\rangle +\sum_{s,t\geq 2}G_{st}\left( z,\overline{z},u\right) 
\hspace{1.2in}\text{for }\alpha =0  \notag \\
v& =-\frac{1}{2\alpha }\ln \left\{ 1-2\alpha \langle z,z\rangle \right\}
+\sum_{s,t\geq 2}G_{st}\left( z,\overline{z},u\right) \quad \text{for }%
\alpha \neq 0  \tag*{(3.1)}  \label{alpha}
\end{align}
where 
\begin{equation*}
\Delta G_{22}\left( z,\overline{z},u\right) =\Delta ^{2}G_{23}\left( z,%
\overline{z},u\right) =0. 
\end{equation*}
Let $\varphi $ be a biholomorphic mapping leaving the $u$-curve invariant
and preserving the form \ref{alpha}. Then $\varphi $ is necessarily given by
the following mapping(cf. \cite{Pa3}): 
\begin{equation*}
\varphi :\left\{ 
\begin{array}{l}
z^{*}=\sqrt{sign\left\{ q^{\prime }(0)\right\} q^{\prime }(w)}Uz\exp \frac{%
\alpha i}{2}\left( q(w)-w\right) \\ 
w^{*}=q\left( w\right)
\end{array}
\right. 
\end{equation*}
where $\alpha \in \Bbb{R}$ and 
\begin{equation*}
\langle Uz,Uz\rangle =sign\left\{ q^{\prime }(0)\right\} \langle z,z\rangle
. 
\end{equation*}
Suppose that the biholomorphic mapping $\varphi $ transforms $M$ to a real
hypersurface $\varphi \left( M\right) $ defined by 
\begin{equation*}
v=-\frac{1}{2\alpha }\ln \left\{ 1-2\alpha \langle z,z\rangle \right\}
+\sum_{s,t\geq 2}G_{st}^{*}\left( z,\overline{z},u\right) 
\end{equation*}
where 
\begin{equation*}
\Delta G_{22}^{*}\left( z,\overline{z},u\right) =\Delta ^{2}G_{23}^{*}\left(
z,\overline{z},u\right) =0. 
\end{equation*}
Then we have the following relation 
\begin{eqnarray*}
q^{\prime }(u)G_{22}^{*}\left( Uz,\overline{Uz},q(u)\right) &=G_{22}\left( z,%
\overline{z},u\right) \\
q^{\prime }(u)\sqrt{\left| q^{\prime }(u)\right| }\exp -\frac{\alpha i}{2}%
\left( q(u)-u\right) G_{23}^{*}(Uz,\overline{Uz},q(u)) &=G_{23}(z,\overline{z%
},u)
\end{eqnarray*}
and 
\begin{eqnarray*}
&\frac{q^{\prime \prime \prime }(u)}{3q^{\prime }(u)}-\frac{1}{2}\left( 
\frac{q^{\prime \prime }(u)}{q^{\prime }(u)}\right) ^{2}+\frac{\alpha ^{2}}{6%
}\left( q^{\prime }(u)^{2}-1\right) \\
&=\frac{1}{6n(n+1)(n+2)}\left\{ q^{\prime }(u)^{2}\Delta
^{3}G_{33}^{*}\left( z,\overline{z},q(u)\right) -\Delta ^{3}G_{33}\left( z,%
\overline{z},u\right) \right\} .
\end{eqnarray*}

Notice that 
\begin{equation}
\frac{q^{\prime \prime \prime }}{3q^{\prime }}-\frac{1}{2}\left( \frac{%
q^{\prime \prime }}{q^{\prime }}\right) ^{2}+\frac{\alpha ^{2}}{6}\left(
q^{\prime 2}-1\right) =0  \tag*{(3.2)}  \label{proj}
\end{equation}
whenever 
\begin{equation*}
\Delta ^{3}G_{33}^{*}\left( z,\overline{z},q\right) dq^{2}=\Delta
^{3}G_{33}\left( z,\overline{z},u\right) du^{2}. 
\end{equation*}
We want to restrict the mapping $\varphi $ so that the function $q(u)$ is a
solution of the ordinary differential equation \ref{proj}. The restriction
on $\varphi $ has to be achieved by requiring an additional condition on the
normal form \ref{alpha}.

We claim that the following choice works 
\begin{align}
\Delta ^{3}G_{33}\left( z,\overline{z},u\right) & =\mathrm{const}.\Delta
^{4}\left( G_{22}\left( z,\overline{z},u\right) \right) ^{2}  \notag \\
\Delta ^{3}G_{33}^{*}\left( z,\overline{z},q\right) & =\mathrm{const}.\Delta
^{4}\left( G_{22}^{*}\left( z,\overline{z},q\right) \right) ^{2}. 
\tag*{(3.3)}  \label{cond33}
\end{align}
Because of the relation 
\begin{equation*}
G_{22}^{*}\left( z,\overline{z},q\right) dq=G_{22}\left( z,\overline{z}%
,u\right) du, 
\end{equation*}
the condition \ref{cond33} gives 
\begin{equation*}
q^{\prime }(u)^{2}\Delta ^{3}G_{33}^{*}\left( z,\overline{z},q(u)\right)
=\Delta ^{3}G_{33}\left( z,\overline{z},u\right) 
\end{equation*}
which yields the ordinary differential equation \ref{proj}.

Hence we define a normal form such that 
\begin{eqnarray*}
v &=\langle z,z\rangle +\sum_{s,t\geq 2}G_{st}\left( z,\overline{z},u\right) 
\hspace{1.2in}\text{for }\alpha =0 \\
v &=-\frac{1}{2\alpha }\ln \left\{ 1-2\alpha \langle z,z\rangle \right\}
+\sum_{s,t\geq 2}G_{st}\left( z,\overline{z},u\right) \quad \text{for }%
\alpha \neq 0
\end{eqnarray*}
where 
\begin{eqnarray*}
\Delta G_{22}\left( z,\overline{z},u\right) &=\Delta ^{2}G_{23}\left( z,%
\overline{z},u\right) =0 \\
\Delta ^{3}G_{33}\left( z,\overline{z},u\right) &=\beta \Delta ^{4}\left(
G_{22}\left( z,\overline{z},u\right) \right) ^{2}\quad \text{for some }\beta
\in \Bbb{R}.
\end{eqnarray*}
We easily see that all normalizations associated to any normal form above
are uniquely determined by some constant initial values.

Chern-Moser normal form is given in the case of $\alpha =\beta =0$ so that 
\begin{equation*}
v=\langle z,z\rangle +\sum_{s,t\geq 2}G_{st}\left( z,\overline{z},u\right) 
\end{equation*}
where 
\begin{equation*}
\Delta G_{22}=\Delta ^{2}G_{23}=\Delta ^{3}G_{33}=0. 
\end{equation*}
Moser-Vitushkin normal form is defined by taking $\alpha \neq 0$ and $\beta
=0$ so that 
\begin{equation*}
v=-\frac{1}{2\alpha }\ln \left\{ 1-2\alpha \langle z,z\rangle \right\}
+\sum_{s,t\geq 2}G_{st}\left( z,\overline{z},u\right) 
\end{equation*}
where 
\begin{equation*}
\Delta G_{22}=\Delta ^{2}G_{23}=\Delta ^{3}G_{33}=0. 
\end{equation*}
We shall see each normal form has its own advantage in applications(cf. \cite
{Pa3}).

\textbf{II}. Burns and Shnider \cite{BS} have reported that the geometric
theory of Chern and Moser \cite{CM74} gives a projective parametrization on
a chain which is different from the parametrization defined by Chern-Moser
normal form(cf. \cite{BFG}). From \ref{4.92}, we obtain 
\begin{align}
& \frac{q^{\prime \prime \prime }}{3q^{\prime }}-\frac{1}{2}\cdot \left( 
\frac{q^{\prime \prime }}{q^{\prime }}\right) ^{2}  \notag \\
& =\frac{1}{6n(n+1)(n+2)}\left\{ q^{\prime 2}\Delta ^{3}H_{33}^{*}(z,%
\overline{z},q(u))-\Delta ^{3}H_{33}(z,\overline{z},u)\right\} . 
\tag*{(3.4)}  \label{4.25}
\end{align}
Note that the solution $q(u)$ in \ref{4.25} is given by 
\begin{equation*}
q(u)=\frac{au+b}{cu+d},\text{ }ad-bc\neq 0,\text{ }a,b,c,d\in \Bbb{R}
\end{equation*}
whenever 
\begin{equation}
\varepsilon \equiv \Delta ^{3}H_{33}(z,\overline{z},u)du^{2}=\Delta
^{3}H_{33}^{*}(z,\overline{z},q)dq^{2}.  \tag*{(3.5)}  \label{4.79}
\end{equation}
The mapping \ref{4.93} effects 
\begin{equation*}
H_{2k}(z,\bar{z},u)=q^{\prime }\left| q^{\prime }\right| ^{\frac{k-2}{2}%
}H_{2k}^{*}(Uz,\overline{Uz},q(u))\quad \text{for }k\geq 2.
\end{equation*}
Thus there are many possible choices satisfying the equalities \ref{4.79}.
To reduce the possible choices, we require that the function $\varepsilon $
is independent of a choice of a chain. Clearly there exists such a function $%
\varepsilon $, for instance, $\varepsilon =0$ so that 
\begin{equation*}
\Delta ^{3}H_{33}^{*}(z,\overline{z},q(u))=\Delta ^{3}H_{33}(z,\overline{z}%
,u)=0.
\end{equation*}
The requirement is also satisfied by a function $\varepsilon $ if we define $%
\varepsilon $ as follows: 
\begin{eqnarray*}
\varepsilon  &\equiv &\Delta ^{3}H_{33}du^{2} \\
\Delta ^{3}H_{33} &=&c\Delta ^{4}(H_{22})^{2},
\end{eqnarray*}
for a constant real number $c\in \Bbb{R}.$ Thus we can define a normal form
similar to Chern-Moser normal form except for replacing the condition $%
\Delta ^{3}H_{33}=0$ with 
\begin{equation}
\sum N_{\alpha \beta \gamma ...}^{\text{ }\alpha \beta \gamma }=\frac{4}{9}%
\sum N_{\alpha \beta ..}^{\text{ }\gamma \delta }N_{\gamma \delta ..}^{\text{
}\alpha \beta }  \tag*{(3.6)}  \label{faran}
\end{equation}
where 
\begin{eqnarray*}
H_{22}(z,\overline{z},u) &=&\sum N_{\alpha \beta \overline{\gamma }\overline{%
\delta }}z^{\alpha }z^{\beta }z^{\overline{\gamma }}z^{\overline{\delta }} \\
H_{23}(z,\overline{z},u) &=&\sum N_{\alpha \beta \gamma \overline{\delta }%
\overline{\rho }\overline{\sigma }}z^{\alpha }z^{\beta }z^{\gamma }z^{%
\overline{\delta }}z^{\overline{\rho }}z^{\overline{\sigma }} \\
\sum N_{\alpha \beta \gamma ...}^{\text{ }\alpha \beta \gamma } &=&\frac{1}{%
(3\cdot 2)^{2}}\Delta ^{3}H_{33} \\
\sum N_{\alpha \beta ..}^{\text{ }\gamma \delta }N_{\gamma \delta ..}^{\text{
}\alpha \beta } &=&\frac{1}{3\cdot 2^{5}}\Delta ^{4}(H_{22})^{2}.
\end{eqnarray*}
Then the condition \ref{faran} gives on a chain the parametrization of the
geometric theory of Chern and Moser(cf. \cite{Fa}).

\end{document}